\documentclass[preprint]{elsarticle}

\usepackage{graphicx}
\usepackage{amssymb,amscd}
\usepackage{bbm}
\usepackage{epsfig}

\newtheorem{thm}{Theorem}[section]
\newtheorem{rem}[thm]{Remark}
\newtheorem{algorithm}[thm]{Algorithm}

\begin{document}

\title{Adjoint method for a tumor invasion PDE-constrained optimization problem using FEM}

\author[famaf,ciem]{A. A. I.~Quiroga \corref{cor1}}
\ead{aiquiroga@famaf.unc.edu.ar}
\cortext[cor1]{Corresponding author: Facultad de Matem\'atica, Astronom\'ia y F\'isica,
  Medina Allende s/n, 5000 C\'ordoba, Argentina}

\author[famaf,ciem]{D. R.~Fern\'andez}

\author[famaf,ciem]{G. A.~Torres}

\author[famaf,ciem]{C. V.~Turner}

\address[famaf]{Facultad de Matem\'atica, Astronom\'ia y F\'isica,
  Medina Allende s/n, 5000 C\'ordoba, Argentina}

\address[ciem]{Centro de Investigaciones y Estudios en Matem\'atica -
  CONICET, Medina Allende s/n, 5000 C\'ordoba, Argentina}

\begin{abstract}
  In this paper we present a method for estimating unknown parameter
  that appear on a non-linear reaction-diffusion model of cancer
  invasion. This model considers that tumor-induced alteration of
  micro-enviromental pH provides a mechanism for cancer invasion. A
  coupled system reaction-diffusion describing this model is given by
  three partial differential equations for the non dimensional spatial
  distribution and temporal evolution of the density of normal tissue,
  the neoplastic tissue growth and the excess concentration of H$^+$
  ions. Each of the model parameters has a corresponding biological
  interpretation, for instance, the growth rate of neoplastic tissue,
  the diffusion coefficient, the reabsorption rate and the destructive
  influence of H$^+$ ions in the healthy tissue.
  
  After solving the forward problem properly, we use the model for the
  estimation of parameters by fitting the numerical solution with real
  data, obtained via in vitro experiments and fluorescence ratio
  imaging microscopy. We define an appropriate functional to compare
  both the real data and the numerical solution using the adjoint
  method for the minimization of this functional.

  We apply Finite Element Method (FEM) to solve both the direct and inverse
  problem, computing the \emph{a posteriori} error.
\end{abstract}


% AMS CLASSIFICATION
  % 35K57 -> reaction-diffusion equation
  % 35Q93 -> PDE-constrained optimization
  % 80M10 -> finite element method

\maketitle

\section{Introduction.}

Cancer is one of the greatest killers in the world although medical
activity has been successful, despite great difficulties, at least for
some pathologies. A great effort of human and economical resources is
devoted, with successful outputs, to cancer research,
\cite{Adam1,Adambellomo,BeChaDe09,BellomoLiMaini,byrne2010dissecting,bellomo2000modelling}.

Some comments on the importance of mathematical modeling in cancer can
be found in the literature. In the work \cite{BellomoLiMaini} the
authors say ``Cancer modelling has, over the years, grown immensely as
one of the challenging topics involving applied mathematicians working
with researchers active in the biological sciences. The motivation is
not only scientific as in the industrial nations cancer has now moved
from seventh to second place in the league table of fatal diseases,
being surpassed only by cardiovascular diseases.''

We use in this work the mathematical analyses first proposed by
\cite{gatenby1996reaction} which supports the acid-mediated invasion
hypothesis, hence it is acquiescent to mathematical representation as a
reaction-diffusion system at the tissue scale, describing the spatial
distribution and temporal development of tumor tissue, normal tissue,
and excess H$^+$ ion concentration.

The model predicts a pH gradient extending from the tumor-host
interface. The effect of biological parameters critical to controlling
this transition is supported by experimental and clinical
observations \cite{martin1994noninvasive}.

In \cite{gatenby1996reaction} the authors model tumor invasion in an
attempt to find a common, underlying mechanism by which primary and
metastatic cancers invade and destroy normal tissues. They are not
modeling the genetic changes which result in transformation nor do
they seek to understand the causes of these changes. Similarly, they
do not attempt to model the large-scale morphological features of
tumors such as central necrosis. Rather, they concentrate on the
microscopic scale population interactions occurring at the tumor-host
interface, reasoning that these processes strongly influence the
clinically significant manifestations of invasive cancer.

Specifically, the authors hypothesize that transformation-induced
reversion of neoplastic tissue to primitive glycolytic metabolic
pathways, with resultant increased acid production and the diffusion
of that acid into surrounding healthy tissue, creates a peritumoral
microenvironment in which tumor cells survive and proliferate, whereas
normal cells are unable to remain viable. The following temporal
sequence would derive: (a) high H$^+$ ion concentrations in tumors
will extend, by chemical diffusion, as a gradient into adjacent normal
tissue, exposing these normal cells to tumor-like interstitial pH; (b)
normal cells immediately adjacent to the tumor edge are unable to
survive in this chronically acidic environment; and (c) the
progressive loss of layers of normal cells at the tumor-host interface
facilitates tumor invasion. Key elements of this tumor invasion
mechanism are low interstitial pH of tumors due to primitive
metabolism and reduced viability of normal tissue in a pH environment
favorable to tumor tissue.

These model equations depend only on a small number of cellular and
subcellular parameters. Analysis of the equations shows that the model
predicts a crossover from a benign tumor to one that is aggressively
invasive as a dimensionless combination of the parameters increases
through a critical value.

\begin{figure}[ht]
\begin{center}
\includegraphics[scale=0.75]{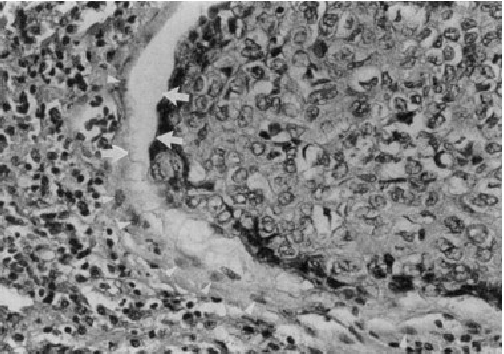} %{tumor_real.eps}
\setlength{\belowcaptionskip}{2pt}
\caption{A micrographs of the tumor-host interface from human
squamous cell carcinomas of the head and neck \cite{gatenby1996reaction}.}
\label{fig:gap}
\end{center}
\end{figure}

The dynamics and structure of the tumor-host interface in invasive
cancers are shown to be controlled by the same biological parameters
which generate the transformation from benign to malignant growth. A
hypocellular interstitial gap, as we can see in Figure \ref{fig:gap} \cite[Figure 4a]{gatenby1996reaction},
at the interface is predicted to occur in some cancers.

In this paper we estimate one of these parameters (the destructive
influence of H$^+$ ions in the healthy tissue) using an inverse
problem. Moreover, via fluorescence ratio imaging microscopy, it is
possible get data about the concentration of hydrogen ions
\cite{martin1994noninvasive}. We propose a framework via a
PDE-constrained optimization problem, following the PDE-based model by
Gatenby \cite{gatenby1996reaction}. In this approach, tumor invasion
is modeled via a coupled nonlinear system of partial differential
equations, which makes the numerical solution procedure quite
challenging.

This kind of problem constitutes a particular application of the
so-called inverse problems, which are being increasingly used in a
broad number of fields in applied sciences. For instance, problems
referred to structured population dynamics \cite{PeZu07}, computerized
tomography and image reconstruction in medical imaging
\cite{van2011source,ZuMa03}, and more specifically tumor growth
\cite{AgBaTu,hogea,knopoff}, among many others.

We solve a minimization problem using a gradient-based method considering
the adjoint method in order to find the derivative of an objective
functional. In this way, we would obtain the best parameter that fits
patient-specific data.

The contents of this paper, which is organized into 9 sections and an
Appendix, are as follows: Section 2 consists in some preliminaries
about the model and the definition of the direct problem. Section 3
deals with the variational formulation of the direct problem.  Section
4 considers the formulation of the minimization problem.  Section 5
introduces the reduced and adjoint problem, deriving the optimality
conditions for the problem. Section 6 finds the derivative of the
solution of a functional with respect to a parameter that does not
appear explicitly in the equation. Section 7 deals with the numerical
solution of the adjoint problem, designing a suitable algorithm to
solve it. In particular, we use the Finite Element Method with a
computation of \emph{a posteriori} error. In Section 8 we show some numerical
simulations to give information on the behavior of the functional and
its dependence on the parameters including the corresponding
tables. Section 9 presents the conclusions and some future work
related to the contents of this paper. In the Appendix we include all
algebraics of Section 6.

Some words about our notation.  We use $\langle\cdot,\cdot\rangle$ to
denote the $L^2$ inner product (the space is always clear from the
context) and we consider the sum of inner products for a cartesian
product of spaces.  For a function $F:V\times U_{ad}\rightarrow
\mathcal{Z}$ such that $(u,\delta_1)\mapsto F(u,\delta_1)$, we denote
by $F'(u,\delta_1)$ the full Fr\'echet-derivative and by
$\frac{\partial F}{\partial u}(u,\delta_1)$ and $\frac{\partial
  F}{\partial \delta_1}(u,\delta_1)$ the partial Fr\'echet-derivatives
of $F$ at $(u,\delta_1)$. For a linear operator $T : V \rightarrow
\mathcal{Z}$ we denote $T^* : \mathcal{Z}^* \rightarrow V^*$ the
adjoint operator of $T$. If $T$ is invertible, we call $T^{-*}$ the
inverse of the adjoint operator $T^*$.

\section{Some preliminaries about the model.}
We present the mathematical model of the tumor-host interface based on
the acid mediation hypothesis of tumor invasion due to
\cite{gatenby1996reaction}. For convenience we reproduce the equations
here, which determine the spatial distribution and temporal evolution
of three fields: $N_1 (x, t)$, the density of normal tissue; $N_2 (x,
t)$, the density of neoplastic tissue; and $L (x, t),$ the excess
concentration of H$^+$ ions. The units of $N_1$ and $N_2$ are
cells/cm$^3$ and excess H$^+$ ion concentration is expressed as a
molarity (M), $x$ and $t$ are the position (in cm) and time (in
seconds), respectively.
\begin{eqnarray}
\frac{ \partial N_1 }{ \partial t } &=& r_1 N_1 \left( 1 - \frac{N_1 }{K_1} \right) - d_1 L  N_1, \label{dim1} \\ 
 \frac{ \partial N_2 }{ \partial t } &=& r_2N_2 \left( 1 - \frac{N_2 }{K_2} \right)+ \nabla \cdot \left( D_{N_2} \left( 1 - \frac{N_1}{K_1} \right) \nabla N_2 \right), \label{dim2} \\ 
 \frac{ \partial L}{ \partial t } &=& r_3N_2 - d_3 L  + D_{N_3} \Delta L, \label{dim3} 
\end{eqnarray}
where the variables are in $\Omega \times [0,T]$. 

In equation (\ref{dim1}) the behavior of the normal tissue is
determined by the logistic growth of $N_1$ with growth rate $r_1$ and
carrying capacity $K_1$, and the interaction of $N_1$ with excess
H$^+$ ions leading to a death rate proportional to $L$. The number
$d_1L$ is the excess acid concentration, dependent death rate in
accord with the well-described decline in the growth rate of normal
cells, due to the reduction of pH from its optimal value of $7.4$. The
constants $r_1$, $d_1$ and $K_1$ have units of $1/$s, l/(M s) and
cells/cm$^3$, respectively.
  
For equation (\ref{dim2}), the neoplastic tissue growth is described
by a reaction-diffusion equation. The reaction term is governed by a
logistic growth of $N_2$ with growth rate $r_2$ and carrying capacity
$K_2$. The diffusion term depends on the absence of healthy tissue
with a diffusion constant $D_{N_2}$. Constants $r_2$, $K_2$ and
$D_{N_2}$ have units of $1/$s, cells/cm$^3$ and cm$^2$/s,
respectively.

In equation (\ref{dim3}), it is assumed that excess H$^+$ ions are
produced at a rate proportional to the neoplastic cell density, and
diffuse chemically. An uptake term is included to take account of the
mechanisms for increasing local pH (e.g., buffering and large-scale
vascular evacuation \cite{gatenby1996reaction}). Constant $r_3$ is the
production rate (M cm$^3$/(cell s)), $d_3$ is the reabsorption rate
(1/s), and $D_{N_3}$ is the H$^+$ ion diffusion constant (cm$^2$/s).

All the parameter values can be found in Table \ref{tabla_parametros}.
\begin{table}[!ht]
  \centering
  \begin{tabular}{cr}
    Parameter & Estimate  \\ \hline
    $K_1$ & $5\times10^7/$cm$^3$ \\
    $K_2$ & $5\times10^7/$cm$^3$ \\
    $r_1$ & $1\times10^{-6}/$s \\
    $r_2$ & $1\times10^{-6}/$s \\
    $D_{N_2}$ & $2\times10^{-10}$cm$^2/$s\\ 
    $D_{N_3}$ & $5\times10^{-6}$cm$^2/$s  \\
    $r_3$ & $2.2\times10^{-17} $M cm$^3/$s \\
    $d_3$ & $1.1\times10^{-4}/$s 
  \end{tabular}
  \caption{Parameter values used in \cite{gatenby1996reaction}.}
  \label{tabla_parametros}
\end{table}

\subsection{Nondimensionalization.}
Following the ideas exposed in \cite{gatenby1996reaction}, and
considering that $\Omega \subset \mathbbm{R}$, the mathematical model
is rescaled and the domain is transformed onto the interval $[0,1]
\times [0,T]$. Hence, let us define the following functions:

\begin{equation} \label{nondimensionalization}
\begin{array}{lcl c lcl}
\displaystyle u_1 &=& \displaystyle \frac{N_1}{K_1} & & \displaystyle u_2 &=& \displaystyle \frac{N_2}{K_2} \\[3mm]
\displaystyle u_3 &=& \displaystyle \frac{L}{L_0} & & \displaystyle \tau &=& \displaystyle r_1 t \\[3mm]
\displaystyle \xi &=& \displaystyle \sqrt{ \frac{r_1}{D_{N_3}} } x & & & &
\end{array}
\end{equation}
where $L_0 = r_3 K_2/d_3$. We will continue denoting $x$ and $t$
instead of $\xi$ and $\tau$, respectively. Using the transformation
(\ref{nondimensionalization}) the dimensionless form of the equations
(\ref{dim1})-(\ref{dim3}) become

\begin{eqnarray}
  \frac{ \partial u_1 }{ \partial t }  &=& u_1( 1 - u_1 ) - \delta_1u_1u_3,  \label{nondim1} \\
  \frac{ \partial u_2 }{ \partial t }  &=& \rho_2u_2( 1 - u_2 ) + \frac{\partial }{\partial x } \left( D_2( 1 - u_1) \frac{\partial u_2}{\partial x }  \right), \label{nondim2} \\
  \frac{ \partial u_3 }{ \partial t }  &=& \delta_3( u_2 - u_3 )  + \frac{\partial^2 u_3 }{\partial x^2 }, \label{nondim3}
\end{eqnarray} 
for $(x,t) \in (0,1) \times (0,T]$, where the four dimensionless
quantities which parameterize the model are given by:
\[
\delta_1 = \displaystyle \frac{d_1 r_3 K_2}{d_3 r_1},\qquad \rho_2 = \displaystyle \frac{r_2}{r_1}, \qquad D_2 = \displaystyle \frac{D_{N_2}}{D_{N_3}}, \qquad \delta_3 = \displaystyle \frac{d_3}{r_1}.
\]

The interaction parameters between different cells (healthy and tumor)
and concentration of H$^+$ are difficult to measure
experimentally. This is the reason for which we propose to estimate
them, so we will focus on $\delta_1$ in this work. The other parameters can be estimated by different techniques (see Table \ref{tabla_parametros}).

\subsection{Initial and boundary conditions.}
At $t = 0$ we will consider the tumor at a certain stage of its
evolution. Hence the initial conditions are:

\begin{eqnarray}
  u_1(x,0) &=& u_1^0(x), \label{init1} \\
  u_2(x,0) &=& u_2^0(x), \label{init2} \\
  u_3(x,0) &=& u_3^0(x), \label{init3}
\end{eqnarray}
for all $x \in [0,1]$.
We assume that the tumor is on the left of the domain, in the sense
that the tumor cells are not moving. Then, for all $t \in [0,T]$, we have

\begin{eqnarray}
  \frac{ \partial u_1 }{ \partial x }(0,t) = 0, & u_1(1,t) = 1, \label{borde1} \\
  \frac{ \partial u_2 }{ \partial x }(0,t) = 0, & u_2(1,t) = 0, \label{borde2} \\
  \frac{ \partial u_3 }{ \partial x }(0,t) = 0, & u_3(1,t) = 0. \label{borde3} 
\end{eqnarray}

From now on, equations (\ref{nondim1})-(\ref{borde3}) will be referred
to as the direct problem.
 
\section{Variational form for the direct problem.}
Using the variational techniques for obtaining the weak solution of
the direct problem
\cite{ladyzhenskaia1988linear,kinderlehrer1987introduction,evans1998partial},
we define the following weak formulation:
\begin{eqnarray}
  0 & = & \int_0^T \int_0^1 \lambda_1 
  \left[
    \frac{ \partial u_1 }{ \partial t } - u_1( 1 - u_1 ) + \delta_1u_1u_3
  \right] dxdt + \label{predebil1} \nonumber \\
   &  & \int_0^T \int_0^1 \lambda_2 
  \left[
    \frac{ \partial u_2 }{ \partial t } - \rho_2u_2( 1 - u_2 ) - \frac{\partial }{\partial x } \left( D_2( 1 - u_1) \frac{\partial u_2 }{\partial x } \right) 
  \right] dxdt+ \nonumber \\
   &  &\int_0^T \int_0^1 \lambda_3
  \left[
    \frac{ \partial u_3 }{ \partial t } - \delta_3( u_2 - u_3 )  - \frac{ \partial^2 u_3 }{\partial x^2}
  \right] dxdt, 
\end{eqnarray}
where $\lambda = (\lambda_1, \lambda_2, \lambda_3)$,
\[
  \lambda_1 ,\lambda_2, \lambda_3 \in W = \left\{ v \in L^2(0,T;H^1_D((0,1))) \mbox{ and } \frac{\partial v }{ \partial t } \in L^2(0,T; (H^1_D((0,1)))^* ) \right\},
\]
$$
L^2(0,T;H^1_D((0,1))) = \left\{ v(x,\cdot) \in L^2((0,T)) \mbox{ and } v(\cdot,t) \in H^1_D((0,1)) \right\}
$$
and 
$$
H^1_D = \left\{ v \in H^1((0,1)) : v = 0 \mbox{ on } \Gamma_D = \{ 1 \}\right\}.
$$

Using integration by parts and boundary condition for
$\lambda$ and $u$ in (\ref{predebil1}) we get the following weak
formulation of (\ref{nondim1})-(\ref{borde3}):
\begin{eqnarray}
  0 & = & \int_0^T \int_0^1 \left( \frac{ \partial u_1 }{ \partial t } \lambda_1 - u_1(1-u_1)\lambda_1 + \delta_1u_1u_3\lambda_1 \right) dxdt + \nonumber\\
   &  & \int_0^T \int_0^1 \left( \frac{ \partial  u_2 }{ \partial t }\lambda_2 - \rho_2u_2(1-u_2)\lambda_2 + D_2( 1 - u_1 ) \frac{ \partial u_2 }{ \partial x } \frac{ \partial \lambda_2 }{ \partial x } \right) dxdt +  \nonumber \\
   &  & \int_0^T \int_0^1 \left( \frac{ \partial u_3 }{ \partial t }\lambda_3 + \delta_3u_3\lambda_3 - \delta_3u_2\lambda_3 + \frac{\partial u_3}{\partial x} \frac{\partial \lambda_3}{\partial x} \right) dxdt .\label{formadebil1} 
\end{eqnarray}

A weak solution $u = [u_1, u_2, u_3]^T \in V = W^3 $ is a function
that satisfies (\ref{formadebil1}) for all $\lambda \in V$ and
$u(x,0)=u^0(x) = [u^0_1(x),u^0_2(x),u^0_3(x)] $.

\section{Formulation of the minimization problem.}
As described above we propose to use an inverse problem technique in
order to estimate $\delta_1$. Function $u$ represents the solution of
the direct problem (the components of $u$ are the state variables of
the problem) for each choice of the parameter $\delta_1$.

Let us assume that experimental information is available during the
time interval $0 \leq t \leq T$. Then, the inverse problem can be
formulated as:
\begin{center}
\textit{Find a parameter $\delta_1$ able to generate data $u = [u_1, u_2, u_3]^T$ that best match the available (experimental) information over time $0 \leq t \leq T$}.
\end{center}

For this purpose, we should construct an objective functional which
gives us a notion of distance between the experimental (real) data and
the solution of the system of PDEs for each choice of the parameter
$\delta_1$.

\begin{figure}
\begin{center}
\includegraphics[scale=0.35]{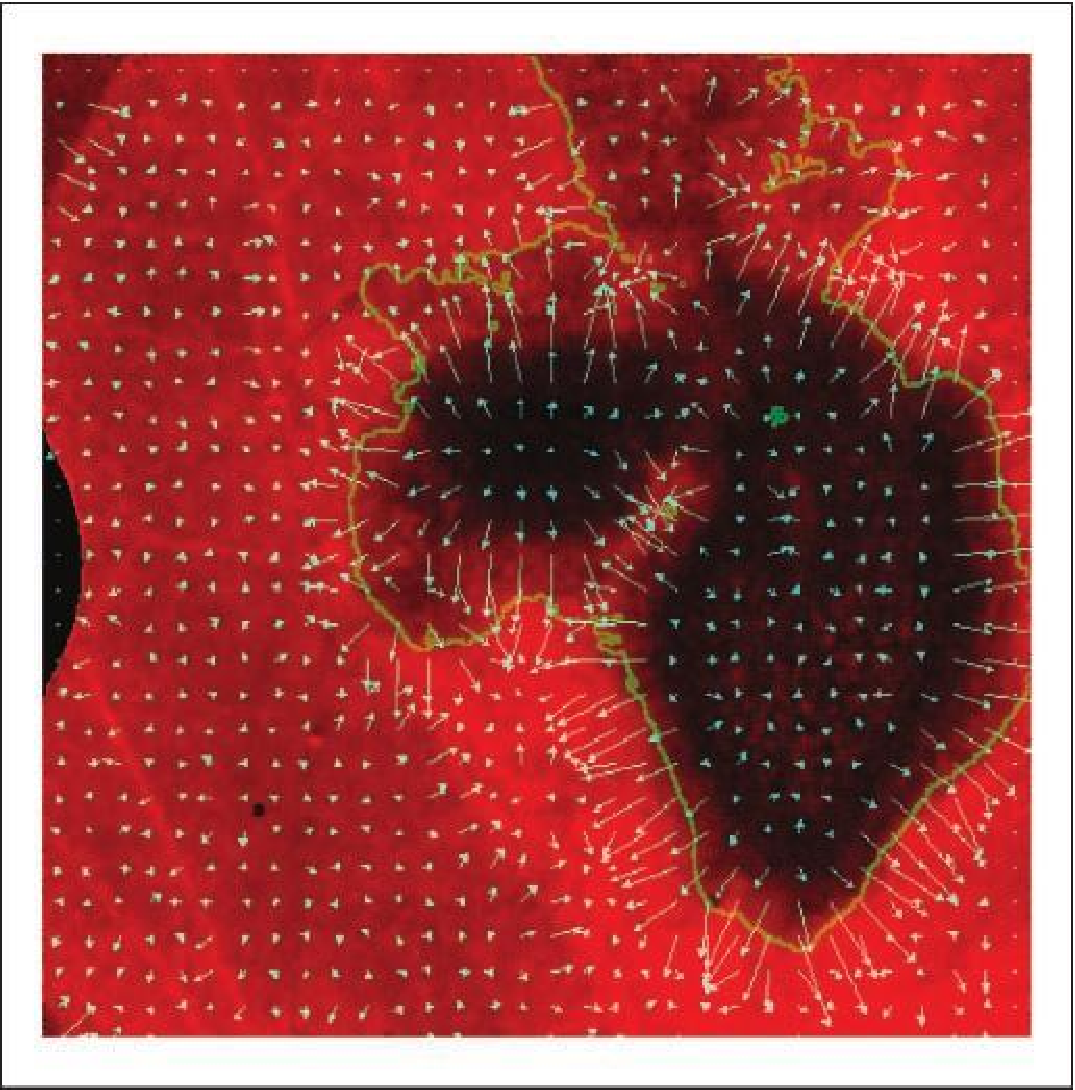} %{acido_real.eps}
\setlength{\belowcaptionskip}{2pt}
\caption{A map of peritumoral H$^+$ flow using vectors generated from
  the pH distribution around the tumor, \cite[Figure 4]{gatenby2006acid}.}
\label{fig:africa}
\end{center}
\end{figure}

First of all, it is important to decide which variables are capable to
be measured experimentally. For instance, the excess concentration of
H$^+$ ions can be measured using fluorescence ratio imaging microscopy
\cite{martin1994noninvasive,gatenby2006acid} at certain times $t_k$,
$k=1,\ldots, M$. For example, Figure \ref{fig:africa} \cite[Figure
4]{gatenby2006acid} shows a map of peritumoral H$^+$ flow using
vectors generated from the pH distribution around the tumor. Such
experiments could help to determine optimal variables and the
parameter in order to control real tumor invasion.

So, the functional $J:V\times U_{ad}\rightarrow \mathbb{R}$ could be defined as:

\begin{equation}
\label{jota1} J(u, \delta_1)=\frac{1}{2} \int_0^T \int_0^1 [u_3(x,t)-\hat{u}_3(x,t)]^2 dxdt,
\end{equation}
where $u_3(x,t)$ is the excess concentration of H$^+$ ions obtained by
solving the direct problem for a certain choice of $\delta_1$ and
$\hat{u}_3(x,t)$ is the excess concentration measured experimentally (real
data).

Let us define $E : V \times U_{ad} \to V^* \times \mathcal{Z}^*$ such that
\begin{eqnarray}  
  \left\langle E(u,\delta_1),\zeta \right\rangle & = & 
  \int_0^T \int_0^1 \left( \frac{ \partial u_1 }{ \partial t } \lambda_1 - u_1(1-u_1)\lambda_1 + \delta_1u_1u_3\lambda_1 \right) dxdt + \nonumber \\ 
  && \int_0^T \int_0^1 \left( \frac{ \partial  u_2 }{ \partial t }\lambda_2 - \rho_2u_2(1-u_2)\lambda_2 + D_2( 1 - u_1 ) \frac{\partial u_2 }{\partial x} \frac{\partial \lambda_2 }{\partial x} \right) dxdt + \nonumber\\
  && \int_0^T \int_0^1 \left( \frac{ \partial u_3 }{ \partial t }\lambda_3 + \delta_3u_3\lambda_3 - \delta_3u_2\lambda_3 + \frac{\partial u_3 }{\partial x} \frac{\partial \lambda_3 }{\partial x} \right) dxdt + \nonumber\\
  &&  \int_0^1 (u_1(x,0) - u^0_1(x))\gamma_1 dx + \int_0^1( u_2(x,0) - u^0_2(x))\gamma_2 dx + \nonumber \\
  && \int_0^1( u_3(x,0) - u^0_3(x))\gamma_3 dx \nonumber \\
  & = & \left\langle \frac{\partial u }{\partial t},\lambda \right\rangle_{V^*,V}  + \left\langle F(u) , \lambda \right\rangle_{V^*,V}
    + \left\langle u(x,0)- u^0(x) , \gamma \right\rangle_{\mathcal{Z}^*,\mathcal{Z}}, \label{constraint}
\end{eqnarray}
where $\zeta = [\lambda,\gamma]$, $\gamma = [ \gamma_1 ,\gamma_2, \gamma_3] \in \mathcal{Z} $ and $\mathcal{Z} = \left(H^1_D((0,1))\right)^3$.

In this way we can rewrite the weak formulation (\ref{formadebil1}) as
$E(u,\delta_1)=0$.

The parameter that best matches the experimental information with the
generated data provided by the direct problem can be computed by solving
a PDE-constrained optimization problem, namely:

\begin{equation}
\begin{array}{rl}
  \displaystyle \mathop{ \mathrm{minimize} }_{\delta_1} & J( u , \delta_1 ) \\ \mathrm{subject \, to} & E(u,\delta_1) = 0, \\ & \delta_1 \in U_{ad},
\end{array}
\end{equation}
where $U_{ad}$ denotes the set of admissible values of $\delta_1$. In
our case, $U_{ad}$ should be a subset of $(0,\infty)$. Notice that a
solution $(u,\delta_1)$ must satisfy the constraint
$E(u,\delta_1)=0$, which constitutes the direct problem.

We remark that, in general, there is a fundamental difference between
the direct and the inverse problems. In fact, the latter is usually
ill-posed in the sense of existence, uniqueness and stability of the
solution. This inconvenient is often treated by using some
regularization techniques
\cite{van2011source,engl1996regularization,kirsch2011introduction}.

\section{Formulation of the reduced and adjoint problems.}
In the following, we will consider the so-called reduced problem

\begin{equation}
  \label{jotamoño}
  \begin{array}{rl}
    \displaystyle \mathop{ \mathrm{minimize} }_{\delta_1} & \tilde{J}(\delta_1)= J(u(\delta_1),\delta_1) \\ \mathrm{subject \, to} & \delta_1 \in U_{ad},
  \end{array}
\end{equation}
where $u(\delta_1)$ is given as the solution of
$E(u(\delta_1),\delta_1)=0$. The existence of the function $u$ is
obtained by the implicit function theorem. According to the ideas
exposed in \cite{brandenburg2009continuous,hinze2009optimization},
this can be done since $U_{ad} = [0,L]$ is a nonempty, closed and
convex set, $J$ and $E$ are continuously Fr\'echet-differentiable
functions, and assuming that for each $\delta_1\in U_{ad}$ there
exists a unique corresponding solution $u(\delta_1)$ such that $E(
u(\delta_1) , \delta_1 ) = 0$ and the derivative $\frac{\partial
  E}{\partial u} (u(\delta_1),\delta_1)$ is a continuous linear
operator continuously invertible for all $\delta_1\in U_{ad}$.

In order to find a minimum of the continuosly differentiable function
$\tilde{J}$, it will be important to compute the derivative of this
reduced objective function. Hence, we will show a procedure to obtain
$\tilde{J}'$ by using the adjoint approach. Since

\begin{equation} \label{jtildedep}
  \tilde{J}^{\ \prime}(\delta_1) = 
  \bigl( u^{\prime}(\delta_1) \bigl)^* \frac{\partial J}{\partial u} (u(\delta_1),\delta_1) + 
  \frac{\partial J}{\partial \delta_1}(u(\delta_1),\delta_1).
\end{equation}

Let us consider $\zeta \in V \times \mathcal{Z}$ as the solution of
the so-called adjoint problem:

\begin{equation}\label{adjoint}
  \frac{\partial J}{\partial u}(u(\delta_1),\delta_1)+\left(\frac{\partial E}{\partial u}(u(\delta_1),\delta_1)\right)^*\zeta =0.
\end{equation}
where $\left( \frac{\partial E}{\partial u}(u,\delta_1) \right)^*$ is
the adjoint operator of $\frac{\partial E}{\partial u}(u,\delta_1)$.
Note that each term in (\ref{adjoint}) is an element of the space
$V^*$.

An equation for the derivative $u^{\prime}(\delta_1)$ is obtained by
differentiating the equation $E(u(\delta_1),\delta_1)=0$ with respect
to $\delta_1$:
\begin{equation}
  \label{derivE}
  \frac{\partial E}{\partial u}(u(\delta_1),\delta_1)u^{\prime}(\delta_1)+\frac{\partial
    E}{\partial \delta_1}(u(\delta_1),\delta_1)=0, 
\end{equation}
where $0$ is the zero vector in $V^* \times \mathcal{Z}^* $.

By using (\ref{jtildedep}) we have that:
\begin{eqnarray*}
  \tilde{J}^{\ \prime}(\delta_1) &=& \bigl( u^{\prime}(\delta_1) \bigl)^* \frac{\partial J}{\partial u} (u(\delta_1),\delta_1) + \frac{\partial J}{\partial \delta_1}(u(\delta_1),\delta_1)
  \\
  & = & - \left( \frac{\partial E}{\partial \delta_1}(u(\delta_1),\delta_1) \right)^* \left( \frac{\partial E}{\partial u}(u(\delta_1),\delta_1) \right)^{-*} \frac{\partial J}{\partial u}(u(\delta_1),\delta_1) + \frac{\partial J}{\partial \delta_1}(u(\delta_1),\delta_1)
  \\
  & = & \left( \frac{\partial E}{\partial \delta_1}(u(\delta_1),\delta_1) \right)^* \zeta + \frac{\partial J}{\partial \delta_1}(u(\delta_1),\delta_1),
\end{eqnarray*}
where in the second equation we used (\ref{derivE}) and for the last
equation we used (\ref{adjoint}). Then:
\begin{equation} \label{adjoint2} \tilde{J}^{\ \prime}(\delta_1) =
  \frac{\partial J}{\partial \delta_1}(u(\delta_1),\delta_1) + \left(
    \frac{\partial E}{\partial \delta_1}(u(\delta_1),\delta_1)
  \right)^* \zeta.
\end{equation}

Notice that in order to obtain $\tilde{J}^{\ \prime}(\delta_1)$ we
need first to compute $u(\delta_1)$ by solving the direct problem,
followed by the calculation of $\zeta$ by solving the adjoint
problem. For computing the second term of (\ref{adjoint2}) it is not
necessary to obtain the adjoint of $\frac{\partial E}{\partial
  \delta_1}(u(\delta_1),\delta_1)$ but just its action over $\zeta$.

\section{Getting the derivative of the functional.}
\label{sec:adjoint}
In order to obtain the adjoint operator of $\frac{\partial E}{\partial u}$, we have to find $\left(\frac{\partial E}{\partial u}\right)^*$ such that:
\begin{equation}
  \label{adjoint3}
  \left\langle \frac{ \partial E }{ \partial u } \eta , \zeta\right\rangle = 
  \left\langle \eta , \left( \frac{ \partial E }{ \partial u } \right)^* \zeta \right\rangle,
\end{equation}
where $ \eta = [\eta_1 , \eta_2, \eta_3]^T$ is the direction of
descent for the state variables $ u_1 $, $ u_2 $ and $ u_3 $, respectively, then
\[
\left\langle \frac{\partial E}{\partial u }( u , \delta_1 ) \eta , \zeta \right\rangle = 
\lim_{\mu\rightarrow 0^+}\frac{ \left\langle E( u + \mu \eta, \delta_1),\zeta \right\rangle - \left\langle E( u , \delta_1 ) , \zeta \right\rangle}{ \mu }.
\]
After some algebraics, it can be shown that $\frac{\partial
  E}{\partial u}\left( u , \delta_1 \right) \eta$ is given by:
\begin{eqnarray}
  \left\langle \frac{\partial E}{\partial u }( u , \delta_1 )\eta , \zeta \right\rangle & = &
  \int_0^T\int_0^1 \left(\frac{ \partial \eta_1 }{ \partial t } - \eta_1( 1 -  2u_1) + \delta_1\eta_1u_3 + \delta_1u_1\eta_3\right)\lambda_1 dxdt + \nonumber \\
  & &\int_0^T\int_0^1 \left(\frac{ \partial \eta_2 }{ \partial t } - \rho_2\eta_2( 1 - 2u_2 )\right)\lambda_2 dxdt +\nonumber \\
  && \int_0^T\int_0^1\left(- D_2 \eta_1 \frac{\partial u_2 }{\partial x} + D_2( 1 - u_1) \frac{\partial \eta_2}{\partial x}\right)\frac{\partial \lambda_2}{\partial x}  dxdt +  \nonumber\\
  & & \int_0^T\int_0^1 \left(\frac{ \partial \eta_3 }{ \partial t } - \delta_3( \eta_2 - \eta_3 )\right)\lambda_3 dxdt + \int_0^T\int_0^1  \frac{\partial \eta_3 }{\partial x}\frac{\partial \lambda_3 }{\partial x}  dxdt +  \nonumber\\
  & &\int_0^1 \eta_1(x,0) \gamma_1(x) dx + \int_0^1 \eta_2(x,0) \gamma_2(x) dx +\int_0^1 \eta_3(x,0) \gamma_3(x) dx.    \label{constraintderiv}  
\end{eqnarray}

An inspection over equations (\ref{adjoint3}) and
(\ref{constraintderiv}) shows that, roughly speaking, we should remove
the spatial and temporal derivatives from $\eta$ and \textit{pass}
them to $\lambda$.

The calculations make use of successive integration by parts to
express each derivative of $\eta$ in terms of a derivative of
$\lambda$. Omitting here the details, that are shown in the Appendix,
we obtain the following expression of the adjoint problem
(\ref{adjoint}), which consists in finding $\lambda \in V $ satisfying
\begin{eqnarray}
  0 & = & \int_0^T\int_0^1 \left( - \frac{ \partial \lambda_1 }{ \partial t }\eta_1 - \eta_1( 1 -  2u_1)\lambda_1 + \delta_1\eta_1u_3\lambda_1 - D_2 \eta_1 \frac{\partial u_2 }{\partial x}  \frac{\partial \lambda_2 }{\partial x} \right) dxdt + \nonumber\\
  & &\int_0^T\int_0^1 \left( - \frac{ \partial \lambda_2 }{ \partial t }\eta_2 - \rho_2\eta_2( 1 - 2u_2 )\lambda_2  + D_2( 1 - u_1) \frac{\partial \lambda_2 }{\partial x}  \frac{\partial \eta_2 }{\partial x} - \delta_3\eta_2\lambda_3 \right) dxdt + \nonumber \\
  & &\int_0^T\int_0^1 \left( - \frac{ \partial \lambda_3 }{ \partial t }\eta_3 + \delta_3\eta_3\lambda_3 + \frac{\partial \lambda_3 }{\partial x} \frac{\partial \eta_3 }{\partial x}  + \delta_1u_1\eta_3\lambda_1 \right) dxdt + \int_0^T\int_0^1\eta_3( u_3 - \hat{u}_3) dxdt  \nonumber \\
  & = & \left\langle - \frac{\partial \lambda}{ \partial t }, \eta \right\rangle_{V^* , V } + \left\langle H(\lambda),\eta \right\rangle_{V^* , V }, \label{adj1}
\end{eqnarray}
for all $ \eta \in V $ and $\lambda(x,T) = 0 $.
As we show in the Appendix we can define $ \gamma(x) = \lambda(x,0) $.

Equation (\ref{adj1}) shall be solved in order to get $\lambda$.
Notice that the adjoint equations are posed backwards in time, with a
final condition at $t=T$, while the state equations are posed
forward in time, with an initial condition at $t=0$.

In order to obtain the derivative of the functional, according to
(\ref{adjoint2}), we must compute the derivative of $E$ with respect
to $\delta_1$.  Since
\[
\left\langle \frac{\partial E}{\partial \delta_1 }( u , \delta_1 ) q , \zeta \right\rangle = 
\lim_{\mu\rightarrow 0^+}\frac{ \left\langle E( u , \delta_1 + \mu q),\zeta \right\rangle - \left\langle E( u , \delta_1 ) , \zeta \right\rangle}{ \mu },
\]
for $q \in U_{ad}$, then
\[
\left\langle \frac{\partial E}{\partial \delta_1 }( u , \delta_1 ) q , \zeta \right\rangle = q\int_0^T \int_0^1 u_1u_3\lambda_1 dxdt.
\]
Thus, since $\frac{\partial J}{\partial \delta_1}=0$, we obtain an
expression for (\ref{adjoint2}), that is
\begin{equation}
  \label{derivadaJ}
  \tilde{J}'(\delta_1)= \left(\frac{\partial E}{\partial \delta_1}(u(\delta_1),\delta_1)\right)^*\zeta = \int_0^T \int_0^1 u_1u_3\lambda_1 dxdt.
\end{equation}

\section{Designing an algorithm to solve the minimization problem.}
\begin{figure}[h!]
  \centering
  \includegraphics[scale=0.43]{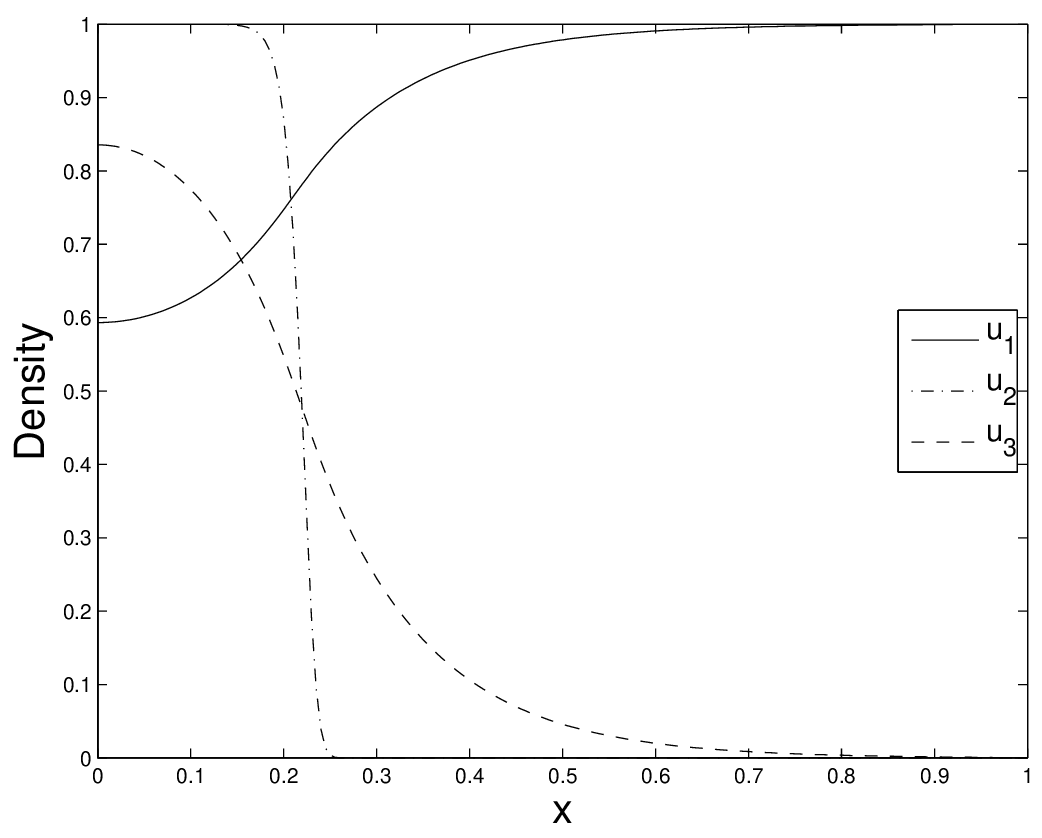}
  \includegraphics[scale=0.43]{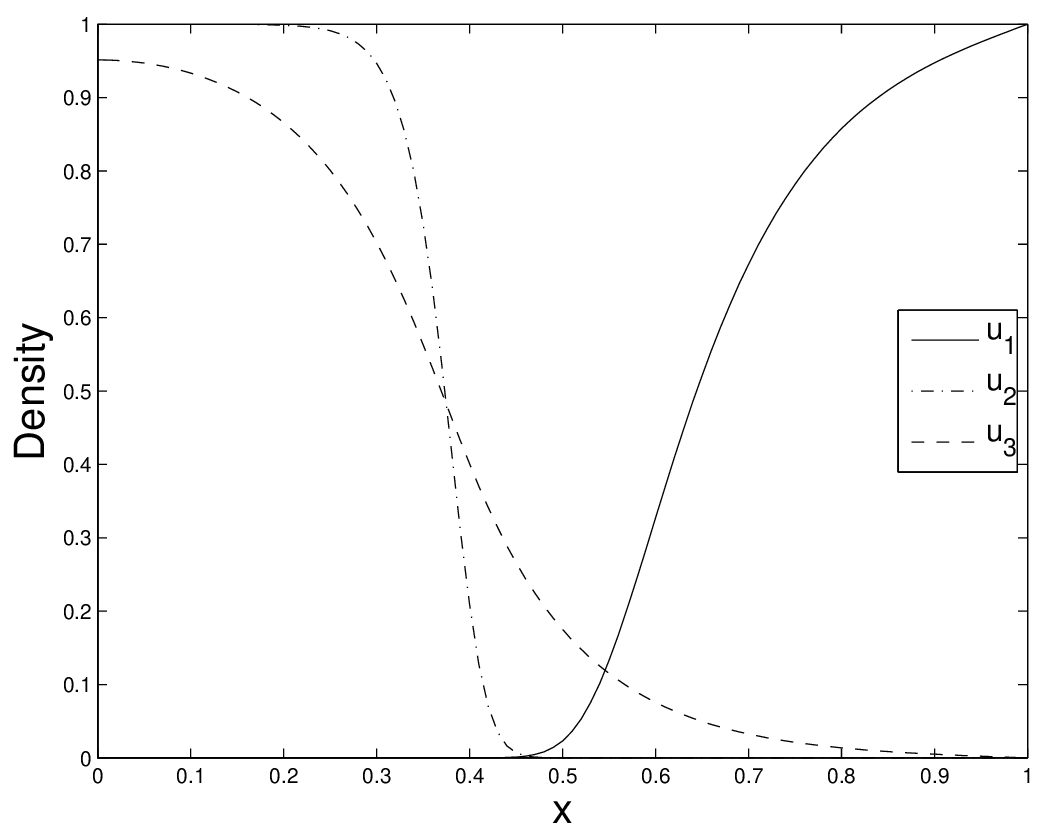}
  \caption{Density of health, tumor cells and excess concetration of H$^+$ at fixed time ($t=20$) in terms of $x$ variable, for $\delta_1 = 0.5$ (left) and $\delta_1 = 12.5$ (right).}
  \label{fig:directS_1}
\end{figure}

It is worth stressing that obtaining model parameters via minimization
of the objective functional $\tilde{J}$ is in general an iterative
process requiring the value of the derivative. To compute $\tilde{J}'$
we just solve two weak PDEs problems per iteration: the direct and the
adjoint problems. This method is much cheaper than the sensitivity
approach \cite{hinze2009optimization} in which the direct problem is
solved many times per iteration. We develop an implementation in
MATLAB that solves the direct and adjoint problems by using a Finite
Element Method and the optimization problem is solved by using a
Sequential Quadratic Programming (SQP) method , using the built-in
function {\tt fmincon}.  For the direct problem, Figure
\ref{fig:directS_1} shows the density of health cells, tumor cells and
excess concetration of H$^+$ at fixed time ($t=20$) in terms of $x$ variable.

It is well-known \cite{nocedal2006numerical} that gradient-based
optimization algorithms require the evaluation of the gradient of the
functional. One important advantage of evaluating the gradient through
adjoints is that it requires to solve the adjoint problem only once
per iteration, regardless the number of inversion variables. Note that
the derivative of the functional can be approximated by using Finite
Element Method.

The method we will use for minimizing the functional $\tilde{J}$ can be summarized as follows:

\begin{algorithm}\label{alg1}
  \textbf{Adjoint-based minimization method.}
  \begin{enumerate}
  \item[1.] Give an initial guess $\delta_1^0$ for the parameter.
  \item[2.] Given $\delta_1^k$ in step $k$, solve the direct and adjoint problems at this step.
  \item[3.] Obtain the derivative of the functional, i.e. $\tilde J'(\delta_1^k)$, using (\ref{derivadaJ}).
  \item[4.] Obtain $\delta_1^{k+1}$ by performing one iteration of the SQP method.
  \item[5.] Stop using the criteria of {\tt fmincon}.
  \end{enumerate}
\end{algorithm}

\begin{algorithm}\label{algdirect}
  \textbf{Direct problem.}
  \begin{enumerate}
  \item[1.] Do an implicit Euler step to find the state variables $u$: $\displaystyle \frac{\partial u}{\partial t}(\cdot,t_n) \approx \frac{ u(\cdot,t_n) - u(\cdot,t_{n-1})}{\tau} = F(u(\cdot,t_n))$, where $t_n = t_{n-1} + \tau $, $F(u(\cdot,t_n))$ is a nonlinear functional and the intial condition is $u^0(x) = u(x,0)$.
  \item[2.] Use FEM  to make a discretization of $u_i(x,t_n)\approx \sum\limits_{j=1}^{nod} u^n_{i,j} \phi_j(x)$, $ i = 1,2,3$, $\phi_j(x)$ are the linear shape function and we note $U_i^n=[u^n_{i,1},\cdots,u^n_{i,j},\cdots,u^n_{i,nod}] \in \mathbb{R}^{nod}$, $U^n=[U_1^n,U_2^n,U_3^n] \in \mathbb{R}^q$, where $ nod $ is the number of uniform distributed nodes for the spatial meshgrid for $[0,1]$.
  \item[3.] Use the Newton method to solve: find $ U^n \in \mathbb{R}^{q} $ such as $U^n - U^{n-1} - \tau G(U^n) = 0 $, where $G$ is the discretization of $F$. 
  \end{enumerate}
\end{algorithm} 

\begin{algorithm}\label{algadj}
  \textbf{Adjoint problem.}
  \begin{enumerate}
  \item[1.] Do an implicit Euler step to find the adjoint variable $\lambda$: $\displaystyle - \frac{\partial \lambda}{\partial t}(\cdot,t_n) \approx - \frac{ \lambda(\cdot,t_n) - \lambda(\cdot,t_{n-1})}{\tau} = H(\lambda(\cdot,t_{n-1}))$, and the final condition is $\lambda(\cdot,T) = 0$.
  \item[2.] Use FEM to make a discretization of $\lambda(\cdot,t_n)$ and solve the linear problem $\lambda^{n-1} - \lambda^n - \tau K(\lambda^{n-1})=0$, where $K$ is the discretization of $H$.
  \end{enumerate}
\end{algorithm} 

\section{Numerical experiments.}

The goal of this section is to test and evaluate the performance of an
adjoint-based optimization method, by executing some numerical
simulations of Algorithm \ref{alg1} for some test-cases.
 
First consider an optimization problem that consist in minimizing the
functional defined in (\ref{jotamoño}), where $\hat{u}_3(x,t)$ is
generated via the forward model, for a choice of the model parameters
$\rho_2 = 1$, $D_2= 4\times 10^{-5}$, $\delta_3 = 1$ and $\hat{\delta}_1 =
0.5, \, 4,\, 12.5 ,\, 16 $. We chose different values of $\hat{\delta}_1$
because each one of these shows a different behavior of tumor
invasion, according to \cite{gatenby1996reaction}.

\begin{figure}[ht]
  \centering
  \includegraphics[scale=.45]{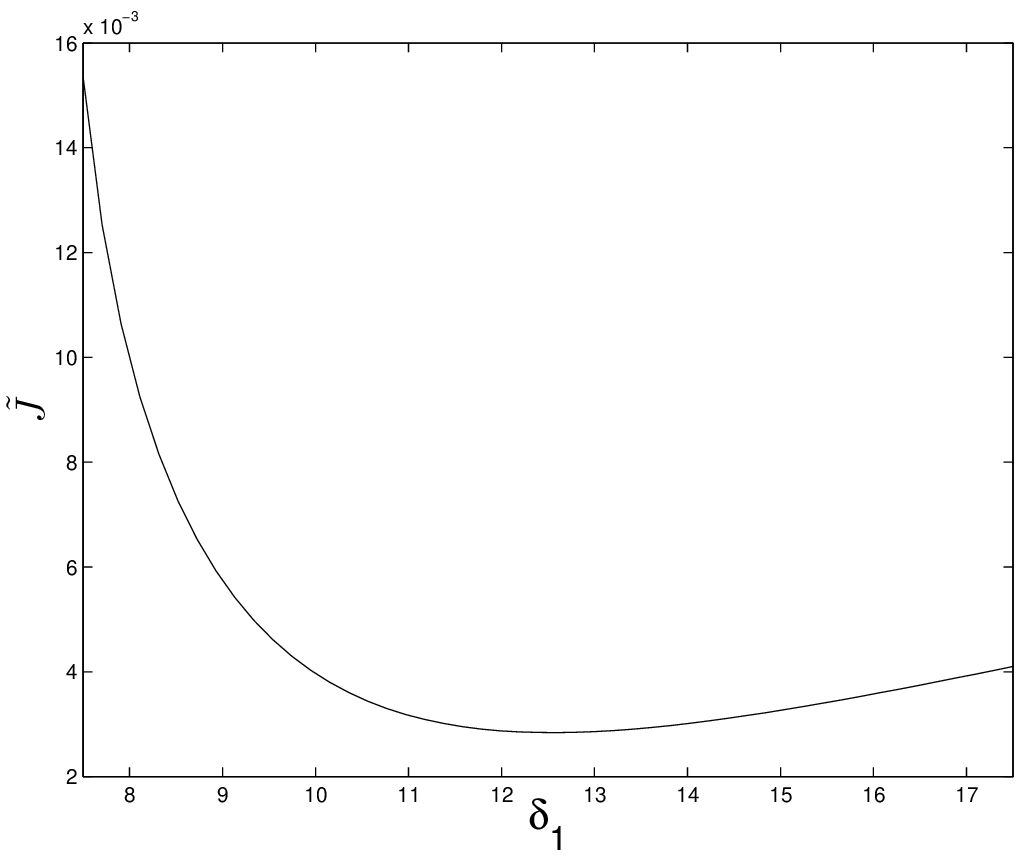}
  \caption{The functional $\tilde{J}$ for $\hat{u}_3$ generated with $\hat{\delta}_1=12.5$.}
  \label{fig:figurej}
\end{figure}

Figure \ref{fig:figurej} shows the value that the functional defined in
(\ref{jotamoño}) takes for different values of $\delta_1$, remaining the
other parameters constant. It is worth mentioning that $\tilde{J}$ looks
convex with respect to $\delta_1$.

The idea of this test case is to investigate how close the original
value of the parameter can be retrieved. However, it is not a
trivial one, because we do not know, for instance, if the
optimization problem has a solution or, in that case, if it is
unique or if the method converges to another local minima.

We have run the Algorithm \ref{alg1} for different values of
$\hat{\delta}_1$ taking the initial condition $\delta_1^0$ randomly,
as we can see from the Table \ref{tabla1} the retrieved parameter
is obtained very accurately since the standard deviation is small.
For Algorithms \ref{algdirect} and \ref{algadj} we use the following
algorithmic parameters $\tau = 0.5 $ and $ T = 20 $, $ nod = 201 $ 
and $U_{ad}=[0,20]$.

\begin{table}[ht]
  \centering
  \begin{tabular}{ r | c | c }
    $\hat{\delta}_1$ & $\bar\delta_1$ & $ S $    \\ \hline
    0.5  & 0.5000 &$\pm$    4.1372$\times10^{-7}$\\
    4    & 4.0000 &$\pm$    2.2187$\times10^{-6}$\\
    12.5 & 12.4999 &$\pm$    4.6521$\times10^{-5}$\\
    16   & 15.9993 &$\pm$    9.4495$\times10^{-5}$
  \end{tabular}
  \caption{Experiments for randomly initial data $\delta_1^0$}
  \label{tabla1}
\end{table}

We emphasize that we have retrieved accurately the value of $\hat{\delta}_1$
independently of the value of $\delta_1^0$. Thus, in the next
experiment we will consider a fixed value $\delta_1^0=8$.

It is well-known that the presence of noise in the data may imply
the appearance of strong numerical instabilities in the solution of
an inverse problem \cite{bertero2006inverse}.

One of the experimental method to obtain values of $\hat{u}_3$ is by using
fluorescence ratio imaging microscopy \cite{martin1994noninvasive}. As
it is well-known that measurements are often affected by perturbations,
usually random ones.

Then we perform numerical experiments where $\hat{u}_3$ is perturbed
by using Gaussian random noise whith mean zero and standard deviation
$\sigma = 0.05,\, 0.1, \, 0.15, \, 0.2 ,\, 0.25, \, 0.3$.  In the next
tables \ref{delta1=05}-\ref{delta1=16}, for each value of $\sigma$, we
show the average $\bar\delta_1$ of 30 values of $\delta_1$, the
standard deviation $S$ and the relative error $e_{\delta_1}=\frac{ |
  \hat{\delta}_1 - \bar{\delta}_1 |}{ \hat{\delta}_1}$.

\begin{table}[ht]
  \centering
  \begin{tabular}{ c | c | c | c }
    $\sigma$ &  $\bar\delta_1$ & $ S $  & $e_{\delta_1}$  \\ \hline
    0.0500 &   0.4707 &$\pm$   0.1231 &   0.0586 \\
    0.1000 &   0.5090 &$\pm$   0.0335 &   0.0180\\
    0.1500 &   0.4855 &$\pm$   0.0472 &   0.0291\\
    0.2000 &   0.4982 &$\pm$   0.0726 &   0.0037\\
    0.2500 &   0.5112 &$\pm$   0.1022 &   0.0225\\
    0.3000 &   0.5027 &$\pm$   0.0937 &   0.0054\\
  \end{tabular}
  \caption{Experiments for $\hat{\delta}_1 = 0.5$}
  \label{delta1=05}
\end{table}

\begin{table}[h!]
  \centering
  \begin{tabular}{ c | c | c | c }
    $\sigma$ &  $\bar\delta_1$ & $ S $  & $e_{\delta_1}$  \\ \hline
    0.0500 &   4.0221 &$\pm$   0.1129 &   0.0055\\
    0.1000 &   4.0470 &$\pm$   0.1695 &   0.0117\\
    0.1500 &   3.9087 &$\pm$   0.2412 &   0.0228\\
    0.2000 &   3.9459 &$\pm$   0.3524 &   0.0135\\
    0.2500 &   3.8970 &$\pm$   0.4800 &   0.0258\\
    0.3000 &   4.0219 &$\pm$   0.4471 &   0.0055\\
  \end{tabular}
  \caption{Experiments for $\hat{\delta}_1 = 4$}
  \label{delta1=4}
\end{table}

\begin{table}[ht]
  \centering
  \begin{tabular}{ c | c | c | c }
    $\sigma$ &  $\bar\delta_1$ & $ S $  & $e_{\delta_1}$  \\ \hline 
    0.0500&   12.7922 &$\pm$   1.2354 &   0.0234\\
    0.1000&   13.0807 &$\pm$   1.9360 &   0.0465\\
    0.1500&   12.0701 &$\pm$   2.3401 &   0.0344\\
    0.2000&   11.4698 &$\pm$   2.7463 &   0.0824\\
    0.2500&   11.1943 &$\pm$   3.7566 &   0.1044\\
    0.3000&   11.8203 &$\pm$   4.3648 &   0.0544\\
  \end{tabular}
  \caption{Experiments for $\hat{\delta}_1 = 12.5 $}
  \label{delta1=12.5}
\end{table}

\begin{table}[ht]
  \centering
  \begin{tabular}{ c | c | c | c }
    $\sigma$ &  $\bar\delta_1$ & $ S $  & $e_{\delta_1}$  \\ \hline
    0.0500 &  16.4165 &$\pm$   2.0834 &   0.0261\\
    0.1000 &  16.6122 &$\pm$   2.7864 &   0.0383\\
    0.1500 &  14.8108 &$\pm$   3.3098 &   0.0743\\
    0.2000 &  13.7965 &$\pm$   3.8915 &   0.1377\\
    0.2500 &  14.1021 &$\pm$   4.5295 &   0.1186\\
    0.3000 &  13.3095 &$\pm$   4.5152 &   0.1681\\
  \end{tabular}
  \caption{Experiments for $\hat{\delta}_1 = 16$}
  \label{delta1=16}
\end{table}

\begin{rem}
 Since we have used FEM to solve both Algorithms,
\ref{algdirect} and \ref{algadj}, we have computed the \emph{a
  posteriori} error in each case \cite{Verfürth_aposteriori,babuvska1978posteriori}. In Table
\ref{tableerror} we put the estimation of the \emph{a posteriori}
error for Algorithm \ref{algdirect} for each $\hat{\delta}_1$.
\end{rem}

\begin{table}[ht]
  \centering
  \begin{tabular}{ r | c | c | c }
    $\hat{\delta}_1 $ & $u_1$ & $u_2$ & $u_3$ \\ \hline
    0.5  & $ 1.72 \times 10^{-14}$ & $ 2.43 \times 10^{-10}$ & $ 2.12 \times^{-7}$ \\
    4    & $ 1.45 \times 10^{-14}$ & $ 4.19 \times 10^{-10}$ & $ 1.80 \times^{-7}$ \\
    12.5 & $ 9.51 \times 10^{-13}$ & $ 1.12 \times 10^{-9}$ & $ 7.61 \times^{-7}$ \\
    16   & $ 5.57 \times 10^{-13}$ & $ 1.06 \times 10^{-9}$ & $ 7.61 \times^{-7}$ \\
  \end{tabular}
  \caption{\emph{A posteriori} error for Algorithm \ref{algdirect}.}
  \label{tableerror}
\end{table}

\section{Final conclusions and future work.}

A miscellany of new strategies, experimental techniques and
theoretical approaches are emerging in the ongoing battle against
cancer. Nevertheless, as new, ground-breaking discoveries relating to
many and diverse areas of cancer research are made, scientists often
have recourse to mathematical modelling in order to elucidate and
interpret these experimental findings, \cite{Adambellomo,
  BellomoLiMaini,byrne2010dissecting,araujo2004history}, and it became
clear that these models are expected to success if the parameters
involved in the modeling process are known. Or eventually, taking into
account that some biological parameters may be unknown (especially in
vivo), the model can be used to obtain them \cite{AgBaTu,
  van2011source}.

This paper, as already mentioned in Section 1, aims at offering a
mathematical tool for the obtention of phenomenological parameters
which can be identified by inverse estimation, by making suitable
comparisons with experimental data. The inverse problem was stated as
a PDE-constrained optimization problem, which was solved by using the
adjoint method. In addition, the gradient of the proposed functional
is obtained and can be extended, in principle, to any number of
unknown parameters.

We remark that the parameter estimation via PDE-constrained
optimization is a general approach that can be used, for instance, to
consider the effects of nonlinear interaction between the health and
tumor cells \cite{mcgillen2013general}.

As a future work we are interested in the dependence of the $\delta_1$
on time and in the dependence of the diffusivity coefficient of excess
of the H$^+$ concentration $D_{N_3}$ with respect to the space
variable $x$, as in \cite{martin2010tumour}. Also we propose to solve
the problem in two dimensional space, where the importance of using
adaptive FEM will be crucial.

\section*{Acknowledgments.}
We appreciate the courtesy of Claudio Padra, from the \textit{Grupo de
  Mec\'anica Computacional - CNEA Bariloche - Argentina}, who strongly
contributed with information above FEM and \emph{a posteriori} error.

The work of  the authors  was partially supported by grants from
CONICET, SECYT-UNC and PICT-FONCYT.

\newpage
\appendix
\section{Appendix: obtaining the adjoint problem.}

In this section we show the calculations involved in order to obtain
the adjoint equations (\ref{adj1}). As stated in Section
\ref{sec:adjoint}, the adjoint equations constitute a weak formulation
of the adjoint problem, with unknown $\zeta$, given by
(\ref{adjoint}). Here, $(\frac{\partial E}{\partial u})^*\zeta$ is
obtained by using (\ref{adjoint3}). In what follows, we shall obtain
equivalent expressions for each of the six terms of the summation
$\langle \frac{\partial E}{\partial u}\eta,\zeta \rangle$, which are
associated with the six constraints given by $E$ in
(\ref{constraint}).

\[
\left\langle \frac{\partial E}{\partial u }( u , \delta_1 )\eta,\zeta \right\rangle= \lim_{\mu\rightarrow 0^+}\frac{ \left\langle E( u + \mu \eta, \delta_1),\zeta \right\rangle - \left\langle E( u , \delta_1 ),\zeta \right\rangle}{ \mu }.
\]

\begin{eqnarray}
  \left\langle \frac{\partial E}{\partial u }( u , \delta_1 )\eta,\zeta \right\rangle & = &
  \int_0^T\int_0^1 \left( \frac{ \partial \eta_1 }{ \partial t }\lambda_1 - \eta_1( 1 -  2u_1)\lambda_1 + \delta_1\eta_1u_3\lambda_1 + \delta_1u_1\eta_3\lambda_1 \right) dxdt + \nonumber \\
  & & \int_0^T\int_0^1 \left( \frac{ \partial \eta_2 }{ \partial t }\lambda_2 - \rho_2\eta_2( 1 - 2u_2 )\lambda_2\right) dxdt + \nonumber \\
  & & \int_0^T\int_0^1 \left( - D_2 \eta_1 \frac{\partial u_2 }{\partial x} \frac{\partial \lambda_2 }{\partial x}  + D_2( 1 - u_1) \frac{\partial \eta_2 }{\partial x} \frac{\partial \lambda_2 }{\partial x} \right) dxdt +  \nonumber\\
  & & \int_0^T\int_0^1 \left( \frac{ \partial \eta_3 }{ \partial t }\lambda_3 - \delta_3( \eta_2 - \eta_3 )\lambda_3 + \frac{\partial \eta_3 }{\partial x} \frac{\partial \lambda_3 }{\partial x} \right) dxdt +  \nonumber\\
  & &\int_0^1 \eta_1(x,0) \gamma_1 dx + \int_0^1 \eta_2(x,0) \gamma_2 dx +\int_0^1 \eta_3(x,0) \gamma_3 dx ,  \label{constraintderiv_ap}  
\end{eqnarray}
using the integration by parts for time, we obtain 
\begin{eqnarray}
  \left\langle \eta , \left(\frac{\partial E}{\partial u }( u , \delta_1 )\zeta \right)^* \right\rangle & = &
  \int_0^T\int_0^1 \left( - \frac{ \partial \lambda_1 }{ \partial t }\eta_1 - ( 1 -  2u_1)\lambda_1\eta_1 + \delta_1u_3\lambda_1\eta_1 \right) dxdt + \nonumber \\
  & & \int_0^T\int_0^1 \left( - \frac{ \partial \lambda_2 }{ \partial t }\eta_2 - \rho_2( 1 - 2u_2 )\lambda_2\eta_2 - \delta_3\lambda_3\eta_2 \right) dxdt + \nonumber \\
    & & \int_0^T\int_0^1 \left( - D_2 \frac{\partial u_2 }{\partial x} \frac{\partial \lambda_2 }{\partial x}\eta_1   + D_2( 1 - u_1) \frac{\partial \lambda_2 }{\partial x} \frac{\partial \eta_2 }{\partial x} \right) dxdt +  \nonumber\\
    & & \int_0^T\int_0^1 \left( - \frac{ \partial \lambda_3 }{ \partial t }\eta_3 + \delta_3\lambda_3\eta_3 +  \frac{\partial \lambda_3 }{\partial x} \frac{\partial \eta_3 }{\partial x} + \delta_1u_1\lambda_1\eta_3 \right) dxdt + \nonumber \\
  & & \int_0^1 \eta_1(x,0) \left(\gamma_1(x) - \lambda_1(x,0) \right) dx + \int_0^1 \eta_1(x,T)\lambda_1(x,T) dx + \nonumber \\
  & & \int_0^1 \eta_2(x,0) \left(\gamma_2(x) - \lambda_2(x,0) \right) dx + \int_0^1 \eta_2(x,T)\lambda_2(x,T) dx + \nonumber \\
  & & \int_0^1 \eta_3(x,0) \left(\gamma_3(x) - \lambda_3(x,0) \right)dx + \int_0^1 \eta_3(x,T)\lambda_3(x,T) dx, \label{constraintderiv_ap_2}  
\end{eqnarray}
then choosing $ \gamma(x) = \lambda(x,0)$ and $\lambda(x,T) = 0 $ for
all $x \in [0,1]$, we obtain the following expression of $\left(
  \frac{\partial E}{\partial u }( u , \delta_1 ) \zeta \right)^*$:
\begin{eqnarray}
  \left\langle \eta , \left( \frac{\partial E}{\partial u }( u , \delta_1 ) \zeta \right)^* \right\rangle & = &
  \int_0^T\int_0^1 \left( - \frac{ \partial \lambda_1 }{ \partial t }\eta_1 - \eta_1( 1 -  2u_1)\lambda_1 + \delta_1\eta_1u_3\lambda_1 \right) dxdt + \nonumber\\
  & &\int_0^T\int_0^1 \left( - \frac{ \partial \lambda_2 }{ \partial t }\eta_2 - \rho_2\eta_2( 1 - 2u_2 )\lambda_2  - \delta_3\eta_2\lambda_3 \right) dxdt + \nonumber \\
  & & \int_0^T\int_0^1 \left( - D_2 \frac{\partial u_2 }{\partial x} \frac{\partial \lambda_2 }{\partial x}\eta_1   + D_2( 1 - u_1) \frac{\partial \lambda_2 }{\partial x} \frac{\partial \eta_2 }{\partial x} \right) dxdt +  \nonumber\\
  & &\int_0^T\int_0^1 \left( - \frac{ \partial \lambda_3 }{ \partial t }\eta_3 + \delta_3\eta_3\lambda_3 + \frac{\partial \lambda_3 }{\partial x} \frac{\partial \eta_3 }{\partial x} + \delta_1u_1\eta_3\lambda_1 \right)dxdt \label{constraintderiv_ap_3}  
\end{eqnarray}

\end{document}